\documentclass[11pt]{amsart}
\usepackage{times}

\usepackage{cite}
\usepackage{color}
\usepackage{algorithmic}
\usepackage[plain]{algorithm}
\usepackage[pdftex]{graphicx}
\usepackage{amsmath}
\usepackage{amsfonts}
\usepackage{fixltx2e}
\usepackage{stfloats}
\usepackage{subfigure}
\usepackage{url}
\usepackage{rotating}
\usepackage{booktabs}
\usepackage{verbatim}

\newtheorem{theorem}{Theorem}[section]

\newtheorem{lemma}[theorem]{Lemma}

\newtheorem{remark}[theorem]{Remark}

\graphicspath{{images/}}

\begin{document}

\title{The Hildreth's Algorithm with Applications to Soft Constraints for User Interface Layout}
\author{Noreen Jamil}
\address{Noreen Jamil\\Department of Computer Science, The University of Auckland, Auckland, New Zealand}
\email{njam031@aucklanduni.ac.nz}

\author{Xuemei Chen}
\address{Xuemei Chen\\Department of Mathematics, University of Missouri, Columbia, MO, USA}
\email{chenxuem@missouri.edu}

\author{Alex Cloninger}
\address{Alex Cloninger\\Applied mathematics program, Yale University, New Haven, CT, USA}
\email{alexander.cloninger@yale.edu}
\thanks{Alex Cloninger is supported by the National Science Foundation under Award No. DMS-1402254}

\begin{abstract}
The Hildreth's algorithm is a row action method for solving large systems of inequalities.
This algorithm is efficient for problems with sparse matrices, as opposed to direct methods such as Gaussian elimination or QR-factorization.  
We apply the Hildreth's algorithm, as well as a randomized version, along with prioritized selection of the inequalities, to efficiently detect the highest priority feasible subsystem of equations.  We prove convergence results and feasibility criteria for both cyclic and randomized Hildreth's algorithm, as well as a mixed algorithm which uses Hildreth's algorithm for inequalities and Kaczmarz algorithm for equalities.
These prioritized, sparse systems of inequalities commonly appear in constraint-based user interface (UI) layout specifications.  The performance and convergence of these proposed algorithms are evaluated empirically using randomly generated UI layout specifications of various sizes.
The results show that these methods offer improvements in performance over standard methods like Matlab's LINPROG, a well-known efficient linear programming solver, and the recent developed Kaczmarz algorithm with prioritized IIS detection.

\end{abstract}

\maketitle

\section{Introduction}
\subsection{Row action methods}
Linear problems are encountered in a variety of fields such as engineering, mathematics and computer science.
Various numerical methods have been proposed to solve these problems, which can be classified into direct and iterative methods.
Direct methods intend to calculate an exact solution in a finite number of steps, whereas iterative methods start with an initial approximation and produce improved approximations in a theoretically infinite sequence whose limit is the exact solution~\cite{Bhatti:Numerical-Analysis}. 

It is observed that iterative methods are often preferable for sparse systems~\cite{Anita:Numerical-MethodsforScientistandEngineers}.
The advantage is that iterative methods spend minimal processing time on coefficients that are zero.
Direct methods, on the other hand, usually lead to fill-in, i.e.\ coefficients change from an initial zero to a non-zero value during the execution of the algorithm, making the processing slower.
Although there are some techniques to minimize fill-in effects, iterative methods are often faster than direct methods for large and sparse problems~\cite{Michele:Preconditioning}.

This paper adapts some particular iterative methods to solve a system of linear equations and inequalities, with the application to User Interface (UI) layout problem, which is a sparse system .
 These methods are the Kaczmarz algorithm~\cite{Kaczmarz:Lettres}, the orthogonal relaxation method (ORM) \cite{M84}, and the Hildreth's algorithm \cite{hildreth, LC80, ID90}, all of which fall into the category of row action methods.   A row action method is an iterative method that uses only one constraint (row) of the system in each iteration. Therefore it is even more preferable in sparse and high dimensional problems, where memory issues might pose challenges to storing or processing the entire system. The Bregman's method \cite{bregman1967relaxation} is another example of row action methods. Another advantage of row action methods is the low complexity for each iteration.   

The Kaczmarz algorithm solves linear equations by successively projecting onto the hyperplane defined by each equality.  The Hildreth's algorithm and the ORM  are parallels of the Kaczmarz algorithm, but for solving linear inequalities.  Both of these algorithms have deterministic and randomized versions for determining the order of the constraints considered.
All three algorithms will be described in detail in Section \ref{sec:row}.
We consider a natural combination of Hildreth's algorithm and the Kaczmarz algorithm to solve a system of mixed linear equalities and inequalities.  In particular, we use Hildreth when an inequality constraint is encountered, and use Kaczmarz when an equation is processed, see Algorithm 1. 
This mixed algorithm is motivated by the UI layout problem.  

\subsection{The UI layout problem and related work}
Constraints are a suitable mechanism for specifying the relationships among objects.
They are used in the area of logic programming, artificial intelligence and UI specification.
They can be used to describe problems that are difficult to solve, conveniently decoupling the description of the problems from their solution.
Due to this property, constraints are a common way of specifying UI layouts, where the objects are widgets and the relationships between them are spatial relationships such as alignment and proportions.
In addition to the relationships to other widgets, each widget has its own set of constraints describing properties such as minimum, maximum and preferred size.

\begin{figure}[htb]
\centering
\includegraphics[width=0.7\columnwidth]{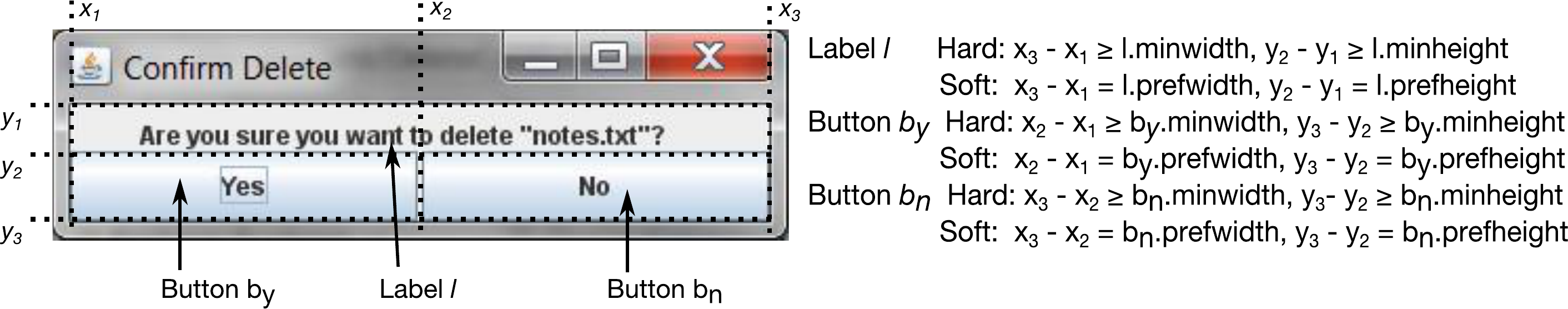}
\caption{Example constraint-based UI layout with hard and soft constraints} \label{fig:Gui-diagram}
\end{figure}

UI layouts are often specified with linear constraints~\cite{Weber:High-Level-Constraints}.
The positions and sizes of the widgets in a layout translate to variables.
Constraints about alignment and proportions translate to linear equations, and constraints about minimum and maximum sizes translate to linear inequalities.
Furthermore, the resulting systems of linear constraints are sparse.
There are constraints for each widget that relate each of its four boundaries to another part of the layout, or specify boundary values for the widget's size, as shown in Figure~\ref{fig:Gui-diagram}.
As a result, the direct interaction between constraints is limited by the topology of a layout, resulting in sparsity.

The linear system in UI layout  problems often contains \textit{conflicting} constraints, i.e. the system is inconsistent.
To resolve conflicts, the notion of \emph{soft constraints} can be introduced.
In contrast to the usual \emph{hard} constraints, which cannot be violated, soft constraints may be violated as much as necessary if no other solution can be found. Many UI layout solvers are based on linear programming and support soft constraints using slack variables in the objective function~\cite{Badros:cassowary,Stuckey:arithmetic-constraints,Finlay:constraint-solving-toolkit,Weber:High-Level-Constraints}.
Most of the direct methods for solving soft constraint problems are least-squares methods such as LU-decompos\-ition and QR-decompo\-sition~\cite{Yasuhiro:least-squares}.
The UI layout solver HiRise~\cite{Hosobe:scalable-linear-constraint} comes under the category of this solver.
HiRise2~\cite{Hosobe:Simplex-Based} is an extended version of the HiRise constraint solver which solves hierarchies of linear constraints by applying an LU-decompos\-ition-based simplex method.

The problem of finding the largest possible subset of constraints that has a feasible solution given a set of linear constraints is widely known as the maximum feasible subsystem (MaxFS) problem~\cite{John:Fast-Heuristics}.
The dual problem to this is the problem of finding the irreducible infeasible subsystem (IIS)~\cite{Amaldi:Maximum-Feasible}.
If one more constraint is removed from an IIS, the subsystem will become feasible.
For both problems, various solving methods are proposed.

There are non-deterministic and deterministic methods to solve the MaxFS problem.
Some of these methods use heuristics~\cite{Amaldi: Two-Phase, O:Feasible-System}, but only a few methods solve the problem deterministically.
The branch and cut method proposed by Pfetsch~\cite{Pfetsch:Branch-Cut} is an example of a deterministic method.

Besides methods for MaxFS there are also some methods to solve the IIS problem.
These methods are: deletion filtering, IIS detection and grouping constraints.
Deletion filtering~\cite{Dravnieks:Minimal} removes constraints from the set of constraints and checks the feasibility of the reduced set.
IIS detection~\cite{Tamiz:IIS} starts with a single constraint and adds constraints successively.
The grouping constraints method~\cite{Guieu:MixedInteger} was introduced to speed up the aforementioned algorithms by adding or removing groups of constraints simultaneously.


In this paper, we prioritize soft constraints  such that, in a conflict between two soft constraints, the soft constraint with the lower priority is violated.
This leads naturally to the notion of \emph{constraint hierarchies}, where all constraints are essentially soft constraints, and the constraints that are considered ``hard'' simply have the highest priorities~\cite{Wilson:Constraint-hierarchies}.
This fits into the UI layout setting well since some constraints can often be more relaxed than others.
Soft constraints guarantee that a problem is always solvable. We naturally modify the IIS detection method to prioritized IIS detection algorithm for solving systems of prioritized linear constraints. In particular, non-conflicting constraints are successively added in descending order of priority until the system becomes infeasible. In each subsystem, we apply the mixed Hildreth's algorithm, both randomly and deterministically, see Algorithm 2. The work  \cite{Gerald:Kaczmarz} also uses prioritized IIS detection method, but with the ORM method. (It is called the Kaczmarz prioritized IIS detection method in \cite{Gerald:Kaczmarz}.) Moreover, they do not consider randomized methods or relaxation parameter. Our numerical experiment shows that our prioritized IIS detection with Hildreth's algorithm converges faster than the one with ORM.

\subsection{Contributions}


The Hildreth's algorithm has never been considered for UI layout problems before. This paper adapts the Hildreth's algorithm to solve systems of linear equations and inequalities, and applies to the UI layout problem. 
Moreover, we provide convergence rate of the algorithm with rigorous analysis, whereas theoretical results about the convergence of LINPROG, LP-Solve and QR-decomposition are still unknown for solving GUI layout problems. To be specific, it has been shown in \cite{ID90} that the cyclic Hildreth's algorithm converges with a linear rate when the system is consistent. 
 This paper extends these results to the mixed Kaczmarz-Hildreth algorithm with a general relaxation parameter. We also consider a randomized version of the mixed Hildreth's algorithm, and  rigorously analyze its convergence rate.


In the case of noise, it has been shown that the method exhibits similar convergence~\cite{Nee10:Randomized-Kaczmarz}, and modified methods even converge to the least-squares solution~\cite{censor1981row,galantai2003projectors,ZF12:Randomized-Extended}.  However, the row action literature has not, to the best of the authors' knowledge, addressed stopping criteria that prevent false convergence for inconsistent systems.  We shall address this point in Section \ref{stopping}.

Our methods are experimentally evaluated with regard to convergence and performance, using randomly generated UI layout specifications.
The results show that both deterministic and randomized mixed Hildreth's algorithm are optimal and efficient in many ways.
We observe that our implementation outperforms Matlab's LINPROG linear optimization package~\cite{linprog:solver}, LP-Solve~\cite{lp:solve}, the implementation of QR-decomposition of the Apache Commons Math Library~\cite{Commons:Math}, and the Kaczmarz prioritized IIS detection method \cite{Gerald:Kaczmarz}.
LP-Solve is a well-known linear programming solver that has been used for UI layout.
The implementation of QR-decomposition of the Apache Commons Math Library is an example of a direct method.  Furthermore, we observe that the randomized version outperforms the deterministic version within the Hildreth's algorithm itself.

The remainder of the paper is organized as follows.
We begin with a background on row action methods for linear systems in Section \ref{sec:row}.
Our algorithms for solving the linear system with prioritized constraints are described in detail in Section \ref{sec:algorithm}.
The convergence result of our algorithm is analyzed in Section \ref{sec:converge}, and  the experimental evaluation are presented in Section~\ref{Evaluation}.


\section{The Kaczmarz algorithm, ORM, and the Hildreth's algorithm }\label{sec:row}
\subsection{Notations}
Throughout this paper, $A$ is an $m\times n$ matrix, whose rows are $a_i$, and $b$ is a vector in $\mathbb{R}^m$. We use $\|A\|_F$ for the Frobenius norm of $A$. Given any vector $z$, $z_i$ is the $i$th component of $z$, and $\|z\|$ is the Euclidean norm of $z$.

We use $Ax=b$ to denote the system of linear equations $\{\langle a_i, x\rangle = b_i\}_{i=1}^m$. With a slight abuse of notation, we use $Ax\leq b$ to denote the system of linear inequalities $\{\langle a_i, x\rangle \leq b_i\}_{i=1}^m$. In the case of a system of linear equations and inequalities (as in the UI problem), we use the notation $Ax=(\leq) b$. No matter which of the three kinds of systems we are in, we assume there is always at least one solution. In the case of $Ax\leq b$ or $Ax=(\leq)b$, the solution set is a convex polytope.  All of these cases are solved via row action algorithms.

Given an index set $J$, let $A_J$ be the submatrix of $A$ that consists of rows indexed by $J$, and $b_J$ be the restriction of $b$ on $J$.

\subsection{The Kaczmarz algorithm}
The Kaczmarz algorithm is a row action method that is used to solve large-scale over-determined linear systems of equations~\cite{Kaczmarz:Lettres}.  It is also used in imaging, tomography, and in that setting is called the ``algebraic reconstruction technique'' (ART)~\cite{GBH70:Algebraic-Reconstruction}.
Given a system of  equations $Ax=b,$
the Kaczmarz algorithm can  be described as follows:
\vspace{0.05in}

\textcolor{white}{lllllllllllllllllllllllllllllllllll}
\begin{tabular}{l}
Initial $x^{(0)}$ (a guess for $x$),\\
$\displaystyle x^{(k+1)}=x^{(k)} +  \alpha\frac{b_{i(k)} - \langle a_{i(k)},  x^{(k)}\rangle }{\|a_{i(k)}\|^2}a_{i(k)}, k\geq0$
\end{tabular}

\vspace{0.05in}

where $\alpha\in (0,2)$ is a relaxation parameter, and  the sequence $i=i(k)$ indicates the order in which the rows are chosen. A common choice is to let $i(k) = (k\mod m) +1$, for which the algorithm cycles through all the rows.
Alternatively, $i$ can be chosen randomly from the discrete set $\{1,2,\cdots, m\}$ according to some probability distribution. For example, the randomized Kaczmarz algorithm~\cite{SV09:Randomized-Kaczmarz} chooses each $i$ independently at random with the distribution $\Pr(i=k)=\|a_i\|^2/\|A\|_F^2$. The algorithm iterates until certain stopping criteria is met, for example, when two consecutive iterates differs less than a pre-specified tolerance.

If $\alpha=1$, at each iteration, the method projects the current iterate $x^{(k)}$ onto the $n-1$ dimensional subspace $\{z: \langle a_i,z\rangle=b_i\}$. If $\alpha>1$, it is called over projection, and if $\alpha<1$, it is called under projection. We observe an accelerated convergence with over projection in our numerical experiments (for Hildreth's algorithm), see Section \ref{Evaluation}.

The convergence rate of the Kaczmarz method depend on the row selection sequence $i(k)$.
A problematic ordering can lead to a drastically reduced rate of convergence.
It is proven in \cite{SV09:Randomized-Kaczmarz} that the randomized Kaczmarz method (when $\alpha=1$) converges to the true solution $x$ in expectation with a linear rate as $\mathbb{E}(\|x^{(k+1)}-x^*\|^2)\leq r \mathbb{E}(\|x^{(k)}-x^*\|^2)$, whre $r<1$ depends on the condition number of $A$.
Random choices of rows (constraints) can overcome a bad ordering of the rows, and is optimal in certain sense \cite{CP2012}.
Randomization is also adopted by ORM or the Hildreth's algorithm as we will see in later sections. A convergence comparison between randomized and deterministic algorithms is discussed in Section \ref{sec:random}.

\subsection{The ORM}
ORM can be considered as a relaxed version of the Kaczmarz method for a system of linear inequalities  $Ax\leq b$. It can be described as

\vspace{0.05in}

\textcolor{white}{lllllllllllllllllllllllllllllllllllllllllllll}
\begin{tabular}{l}
Initial $x^{0}$\\
$\displaystyle c_k = \min\{0, \alpha\frac{b_{i(k)} - \langle a_{i(k)},  x^{(k)}\rangle }{\|a_{i(k)}\|^2}\},$\\
$x^{(k+1)}=x^{(k)} +  c_ka_{i(k)}$
\end{tabular}

\vspace{0.05in}

At each iteration, if $x^{(k)}$ has already satisfied the $i$th constraint, then simply nothing is done; otherwise, an orthogonal projection is performed.

Same as the Kaczmarz algorithm, $i(k)$ can be chosen cyclically or randomly. When chosen cyclically, $\{x^{(k)}\}$ converges to a point $x^*$ on the boundary of $\{x: Ax\leq b\}$, with a linear rate as well~\cite{M84}.
 Leventhal et al~\cite{LL10:Randomized-Methods} generalized the result in~\cite{SV09:Randomized-Kaczmarz}, and proved that the mean square error also converges with a linear rate.
This will be further discussed in Remark \ref{rem:orr}.

\subsection{The Hildreth's algorithm}\label{sec:hildreth}
The Hildreth's algorithm also solves a system of linear inequalities, but with one more benefit: finding the closest point in the solution set to a given point, i.e., it solves the following problem
\begin{align}\label{equ:h1}
&x^*=\arg\min \|x-x^{0}\|\\\label{equ:h2}
&\text{subject to } Ax\leq b.
\end{align}

The algorithm is slightly more complicated than ORM. It is defined as

\vspace{0.05in}

\textcolor{white}{lllllllllllllllllllllllllllllllllllllllllllll}
\begin{tabular}{l}
Initial $x^{(0)}=x^0$, $z^{(0)}=0$,\\
$\displaystyle c_k = \min\{z^{(k)}_{i(k)}, \alpha\frac{b_{i(k)} - \langle a_{i(k)},  x^{(k)}\rangle }{\|a_{i(k)}\|^2}\},$\\
$x^{(k+1)}=x^{(k)} +  c_ka_{i(k)},$ \\
$z^{(k+1)}_{i}=\left\{\begin{array}{lc}z^{(k)}_i, &i\neq i(k)\\z^{(k)}_{i(k)}-c_k, & i=i(k)\end{array}\right.$
\end{tabular}

\vspace{0.05in}

The vector $z^{(k)}$ only gets updated at $i$th component, where $i$ is the active constraint in that iteration. It can be shown that all entries of $z^{(k)}$ are never negative~\cite{LC80}.
The paper \cite{LC80} also provides a geometric interpretation of this algorithm when $\alpha=1$. When the constraint is violated, we project the last approximate to $\{x:\langle a_i,x\rangle= b_i\}$, just like in ORM; when the constraint is satisfied, we under project the last approximate to $\{x:\langle a_i,x\rangle= b_i\}$. The definition of $c_k$ guarantees no over projection.

The sequence $\{x^{(k)}\}$ converges to $x^*$ if appropriate $i(k)$ is chosen, for example, almost cyclic sequences~\cite{LC80}, but no convergence rate is given. Iusem et al~\cite{ID90} established a linear convergence rate with almost cyclic sequences. But so far, no convergence analysis has been given for the randomized version. We shall establish that in Theorem \ref{thm:Hrc}.


\section{Description of the algorithm}\label{sec:algorithm}
In the UI problem, we have a system of linear equations and inequalities $Ax=(\leq)b$, so we  use the mixed Kaczmarz-Hildreth algorithm as mentioned. For convenience, we will call this mixed Hildreth's algorithm Hildreth's algorithm in this and subsequent sections.

\vspace{0.05in}

\begin{tabular}{l}
\hline
\hline
\textbf{Algorithm 1} The Hildreth's algorithm for $Ax=(\leq)b$\\
\hline
1. Initial $x^{(0)}=x^0$, $z^{(0)}=0$,\\
2. Select  $i(k)$,\\
3. if constraint $i(k)$ is an equation, perform Kaczmarz: $\displaystyle x^{(k+1)}=x^{(k)} +  \frac{b_{i(k)} - \langle a_{i(k)},x^{(k)}\rangle }{\|a_{i(k)}\|^2}a_{i(k)}$,\\
\textcolor{white}{3.} if constraint $i(k)$ is an inequality, perform Hildreth:\\
\qquad\qquad$\displaystyle c_k = \min\{z^{(k)}_{i(k)}, \alpha\frac{b_{i(k)} - \langle a_{i(k)},  x^{(k)}\rangle }{\|a_{i(k)}\|^2}\}$\\
\qquad\qquad$x^{(k+1)}=x^{(k)} +  c_ka_{i(k)}$\\
\qquad\qquad$z^{(k+1)}_{i}=\left\{\begin{array}{lc}z^{(k)}_i, &i\neq i(k)\\z^{(k)}_{i(k)}-c_k, & i=i(k)\end{array}\right.$\\
\hline

\end{tabular}
\vspace{0.05in}

During step 2, if each row is chosen cyclically as $i(k)=(k \mod m)+1$, then the resulted method is called \emph{the cyclic Hildreth's algorithm for $Ax=(\leq)b$}.
If $i$ is chosen at random from $\{1,2,\cdots,m\}$ with the distribution $\Pr(i=k)=\|a_i\|^2/\|A\|_F^2$, then the resulted method is called \emph{the randomized Hildreth's algorithm for $Ax=(\leq)b$}.

Both algorithms stop when a string of $N$ consecutive iterations differ by less than a pre-specified tolerance (see Section \ref{stopping}).

As mentioned earlier, we use soft constraints to deal with conflicting constraints.
Some lower priority constraints are allowed to be violated to keep a feasible system.
Every constraint is associated with a priority number, where a bigger priority number indicates a stronger need to be not violated. For example, the priority number for a hard constraint is $\infty$.

We define an index set $E$, which keeps track of \emph{enabled constraints}.
This algorithm starts with an empty set $E$ of enabled constraints.
It then adds one constraint per iteration in order of descending priority to $E$.
In each iteration, the algorithm adds a constraint tentatively to $E$ (``enabling'' it), and attempts to solve the resulting system using Algorithm 1.
If a solution is found, i.e., the subsystem is feasible, the constraint is kept.
Otherwise, the added constraint is removed (index of this constraint is removed from  $E$), in which case the previous solution is restored.
Finally, the algorithm proceeds to the constraint with the next lower priority, until all constraints have been considered. Therefore it is very important to efficiently determine whether a system is feasible/consistent, as we illustrate in Section \ref{stopping}.

\vspace{0.05in}
\textcolor{white}{llllllllllllll}
\begin{tabular}{l}
\hline
\hline
\textbf{Algorithm 2}  Prioritized IIS detection method with Hildreth's algorithm\\
\hline
1. Initiate $E=\emptyset, x^{(0)}=0.$\\
2. for $i=1:m$\\
\quad\quad let $l=l(i)$ be the index of the constraint with $i$th biggest priority number;\\
\quad\quad $E=E\cup\{l\}$;\\
\quad\quad solve $A_Ex=(\leq)b_E$ with Algorithm 1;\\
\quad\quad if $A_Ex=(\leq)b_E$ is feasible (using Algorithm 3) with solution $x$, \\
\quad\quad\quad $x^{(l)}=x;$\\
\quad\quad else\\
\quad\quad\quad $E=E-\{l\}, x^{(l)}=x^{(l-1)}$,\\
\hline
\end{tabular}

\vspace{0.05in}

One can choose either cyclic or randomized Hildreth's algorithm in the ``solving $A_Ex=(\leq)b_E$'' step. The convergence rate of this algorithm will be analyzed in Section \ref{sec:converge}.

\subsection{Feasibility criteria}\label{stopping}
We consider the problem of deciding whether a system is consistent using the randomized version of Algorithm 1. To the best of our knowledge, the literature has never considered this before. First of all, it is easy to show that the system is consistent if and only if Algorithm 1 admits a convergent sequence. However, this is challenging because the inconcsistent system will have convergent subsequences. In fact, for an inconsistent system, it is possible to have $\|x^{(k+i)} - x^{(k+i-1)}\|_{2}<\epsilon$ for $i\in\{1, ..., C\}$ despite $\{x^{(k)}\}$ being a divergent series.  For ease of notation, we shall consider $Ax=b$, though the same argument applies to the case of $Ax=(\le)b$.

With our prioritized IIS detection method, we consider a system $\widetilde{A} = \begin{pmatrix}A \\ a\end{pmatrix}$ and $\widetilde{b}=\begin{pmatrix}\overrightarrow{b}\\ b\end{pmatrix}$, $A\in\mathbb{R}^{m-1 \times n}$ and $\overrightarrow{b}\in\mathbb{R}^{m-1}$, $a\in \mathbb{R}^{n}$ and $b\in \mathbb{R}$.  Let  $Ax= \overrightarrow{b}$ be consistent, but $\widetilde{A}x =\widetilde{b}$ be inconsistent.  Let $x^*$ be the solution to $Ax= \overrightarrow{b}$.

We define $H = \{x: \langle a, x\rangle = b\}$ and $H^* = \{x : \langle a, x\rangle = b + r\}$ where $r = \langle a, x^{*}\rangle - b$.  This makes $H^{*}$ the hyperplane that goes through the solution $x^{*}$.  One can also write $H^{*} = \{w + r\frac{a}{\|a\|} : w\in H\}$.

Let $Q$ be the projection onto $H$, and $Q^{*}$ be the projection onto $H^{*}$.  Then clearly
\begin{eqnarray*}
Qx_{k} - x^{*} = Q^{*}x_{k} - x^{*} - r\frac{a}{\|a\|},
\end{eqnarray*}
which means that by the orthogonality of $a$ to $H$ and $H^{*}$, we have
\begin{eqnarray}\label{eq:projectDistance}
\|Qx_{k} - x^{*}\|^{2} = \|Q^{*}x_{k} - x^{*}\|^{2} + r^{2} \ge r^{2}.
\end{eqnarray}

Equation \eqref{eq:projectDistance} implies that, upon projecting onto the inconsistent equation $ax=b$, the new point $x_{k+1} = Qx_{k}$ will always move at least $r$ away from $x^{*}$, the accumulation point for $Ax = \overrightarrow{b}$.  For this reason, the Hildreth's algorithm will never converge to $x^{*}$.

Detecting these jumps is the key to determining whether the Hildreth's algorithm has converged or is inconsistent.  As we do not have a priori knowledge of $x^{*}$, the best alternative is to examine $\|x_{k+1}-x_{k}\|$.  If the Hildreth's algorithm converged, as we shall see in Theorem \ref{thm:Hrc}, then $\|x_{k+1} - x_{k}\|$ goes to 0 in expectation.

The proposed feasibility criteria is as follows:

\vspace{0.05in}
\textcolor{white}{lllllllllllllllllllllllll}
\begin{tabular}{l}
\hline
\hline
\textbf{Algorithm 3} Feasibility criteria\\
\hline
Input: $maxIterations$, $N$, $\epsilon$,\\
Initial $x^{(0)}=x^0$, $z^{(0)}=0$, counter $=0$,\\
while counter $\le N$ and $k<maxIterations$,\\
\quad Select $i(k)$,\\
\quad Perform randomized Algorithm 1 to generate $x^{(k+1)}$,\\
\quad if $\|x^{(k+1)} - x^{(k)}\|<\epsilon$,\\
\quad\quad counter++, \\
\quad else\\
\quad\quad counter $=0$.\\
\hline
\end{tabular}

\vspace{0.05in}

By choosing an appropriate $N$ and $maxIterations$, we determine whether $\{x^{(k)}\}$ converges or continues to ``jump'' by a distance of $r$ when $i(k)$ corresponds to an inconsistent row $ax=b$, as predicted by \eqref{eq:projectDistance}.  This will determine whether a system is consistent or inconsistent.

For the cyclic version of the algorithm, the choice of $N$ is immediate.  Choosing $N=m$ guarantees $i(k)$ cycles through every row of $A$.  If $\|x^{(k+1)} - x^{(k)}\|<\epsilon$ for each $i\in \{1, ... , m\}$, then the solution $x$ is clearly feasible, and $Ax=(\le) b$ is consistent.

For the randomized version of the algorithm, the choice of $N$ can be found by examining the probability of missing equation $a$ for $N$ consecutive iterations.  This probability is
\begin{eqnarray}\label{eq:prob}
\begin{aligned}
\textnormal{P}(\textnormal{skip row $a$ for $N$ iterations}) &=  \textnormal{P}(\textnormal{skip row $a$})^{N} \\
&= \left(1 - \frac{\|a\|^{2}}{\|A\|_{F}^{2}}\right)^{N}.
\end{aligned}
\end{eqnarray}
One can choose $N$ to match the desired probability of failure.

If the counter reaches $N$ consecutive steps without a jump, then that implies that, with high probability according to \eqref{eq:prob}, there are no inconsistent equations in the system.  Adding this feasibility criteria to the Hildreth algorithm allows us to determine whether the system is consistent.



\section{Convergence result}\label{sec:converge}
This section focuses on the convergence analysis of the Algorithm 1 for solving the following problem:
\begin{align}
&x^*=\arg\min \|x-x^{(0)}\|\\
&\text{subject to } Ax=(\leq) b.
\end{align}

To be more specific, let the system $Ax=(\leq)b$ be $\left\{\begin{array}{cl}\langle a_i, x\rangle \leq b_i & (i\in I_{\leq})\\
\langle a_i, x\rangle =b_i & (i\in I_=)\end{array}\right.$.

Define each hyperplane $H_i=\{x:\langle a_i,x\rangle=b_i\}$, and $I=\{i:x^*\in H_i\}\supset I_=$. Let $S=\bigcap_{i\in I}H_i$ and
$$\mu=\inf_{x\not\in S}\frac{\max_{i\in I}d(x,H_i)}{d(x,S)},$$ where $d(x,X)$ is the distance of $x$ to a set $X$.

We will analyze the convergence rate for both cyclic and randomized versions.
\begin{theorem}
The cyclic Algorithm 1 converges to $x^*$ at a linear rate as
\begin{equation}\label{equ:cHr}
\|x^{(k+m)}-x^*\|^2\leq \frac{1}{1+\sigma}\|x^{(k)}-x^*\|^2,
\end{equation}
where $\sigma=\frac{(2-\alpha)\alpha\mu^2}{1+\alpha^2(m-1)}$.
\end{theorem}
\proof The proof is identical to the proof of Theorem 1 in \cite{ID90} even though \cite{ID90} is for solving \eqref{equ:h1}-\eqref{equ:h2} with the Hildreth's algorithm presented in Section \ref{sec:hildreth}.\qed

\begin{remark}
Notice the rate in \eqref{equ:cHr} is compared after $m$ iterations so that the process goes through all constraints once.
\end{remark}

Before we present the convergence theorem for the randomized Hildreth algorithm, we need a few tools.
\begin{theorem}[Hoffman \cite{hoffman}] Let $S_b$ be the set of feasible solutions of the linear system $Ax=(\leq)b$. Then there exists a constant $L$, independent of $b$, with the following property:
$$\|z-x^*\|\leq L\|e(Az-b)\|,$$
where $e: \mathbb{R}^m\rightarrow\mathbb{R}^m$ is defined by $e(y)=\left\{\begin{array}{cc}
\max\{0,y_i\}& (i\in I_{\leq})\\
y_i & (i\in I_=)\end{array}\right.$.

\end{theorem}
\begin{lemma}[Iusem \cite{ID90}]\label{lem:Herror}
For large enough $k$,
\begin{equation}\label{equ:Herror}
\|x^{(k+1)}-x^*\|^2\leq\|x^{(k)}-x^*\|^2-(\frac{2}{\alpha}-1)\|c_ka_{i(k)}\|^2.
\end{equation}
\end{lemma}

Again even though the results in \cite{ID90} deals with \eqref{equ:h1}-\eqref{equ:h2}, the proof is the same.

\begin{theorem}\label{thm:Hrc}
For large enough $k$, the randomized Algorithm 1  converges to $x^*$ linearly in expectation:
\begin{equation}\label{equ:Hexp}
\mathbb{E}(\|x^{(k+1)}-x^*\|^2)\leq \left(1-\frac{2/\alpha-1}{L^2\|A\|_F^2}\right)\mathbb{E}(\|x^{(k)}-x^*\|^2),
\end{equation}
where $L$ is the Hoffman constant.
\end{theorem}
\proof
At iteration $k$, we use $i=i(k)$. If $x^{(k)}$ violates the $i$th constraint, i.e., $\langle a_i, x^{(k)}\rangle >b_i$, then
$c_ka_i=-\frac{e(Ax^{(k)}-b)_i}{\|a_i\|^2}a_i$.

If $x^{(k)}$ doesn't violate the $i$th constraint, then
$\|c_ka_i\|\geq z_k\|a_i\|\geq0=\left\|\frac{e(Ax^{(k)}-b)_i}{\|a_i\|^2}a_i\right\|$. (It is proved in \cite{LC80} that all entries of $z_k$ are nonnegative.)

Overall, we have $\|c_ka_i\|\geq \left\|\frac{e(Ax^{(k)}-b)_i}{\|a_i\|^2}a_i\right\|$

By Lemma \ref{lem:Herror},
$$\|x^{(k+1)}-x^*\|^2\leq\|x^{(k)}-x^*\|^2-(\frac{2}{\alpha}-1)\frac{e(Ax^{(k)}-b)_i^2}{\|a_i\|^2}.$$
With $x_k$ fixed, and taking expectation over $i(k)$,
\begin{align*}\mathbb{E}(\|x^{(k+1)}-x^*\|^2|x^{(k)})&\leq \|x^{(k)}-x^*\|^2-(\frac{2}{\alpha}-1)\sum_{l=1}^m\frac{\|a_l\|^2}{\|A\|_F^2}\frac{e(Ax^{(k)}-b)_l^2}{\|a_l\|^2}\\
&=\|x^{(k)}-x^*\|^2-(\frac{2}{\alpha}-1)\frac{\|e(Ax^{(k)}-b\|^2}{\|A\|_F^2}\\
&\leq \|x^{(k)}-x^*\|^2-(\frac{2}{\alpha}-1)\frac{\|x^{(k)}-x^*\|^2}{L^2\|A\|_F^2}.
\end{align*}
Taking expectation with respect to $x^{(k)}$ on both sides yields \eqref{equ:Hexp}.
\qed

\begin{remark}\label{rem:orr}
The proof of Theorem \ref{thm:Hrc} heavily relies on the Hoffman bound, and this technique can be found in \cite{LL10:Randomized-Methods}. In fact, the same rate has been shown in \cite{LL10:Randomized-Methods}, but only with $\alpha=1$.
\end{remark}

\subsection{Rate comparison between cyclic and randomized Hildreth's algorithm}\label{sec:random}

We compare the rate between cyclic and randomized Hildreth's algorithm when $\alpha=1$, i.e. we compare these two numbers
$$r_c=\frac{1}{1+\mu^2/m}, r_r=\left(1-\frac{1}{L^2\|A\|_F^2}\right)^m.$$

It is generally difficult to compute either $\mu$ or $L$. However, an estimate can be obtained when the system consists only equations ($I_{\leq}=\emptyset$), in which case the algorithm reduces to the Kaczmarz algorithm. In this situation, one can estimate that $\mu\geq \sigma_{\min}(A)/\|A\|_F$ \cite{Goffin}, and $L=1/\sigma_{\min}(A)$.

In this case, it can be shown that $r_{r}<r_{c}$ as $m\rightarrow\infty$, which shows that the randomized version is more efficient. This is consistent with our numerical experiments.

\section{Experimental Evaluation} \label{Evaluation}

In this section, we present an experimental evaluation of the proposed algorithm.
We conduct two different experiments to evaluate (i) the convergence behavior, and (ii) the performance in terms of computation time.
The experiments are described in the following.

\begin{figure}[t]
\centering
\includegraphics[width=\columnwidth]{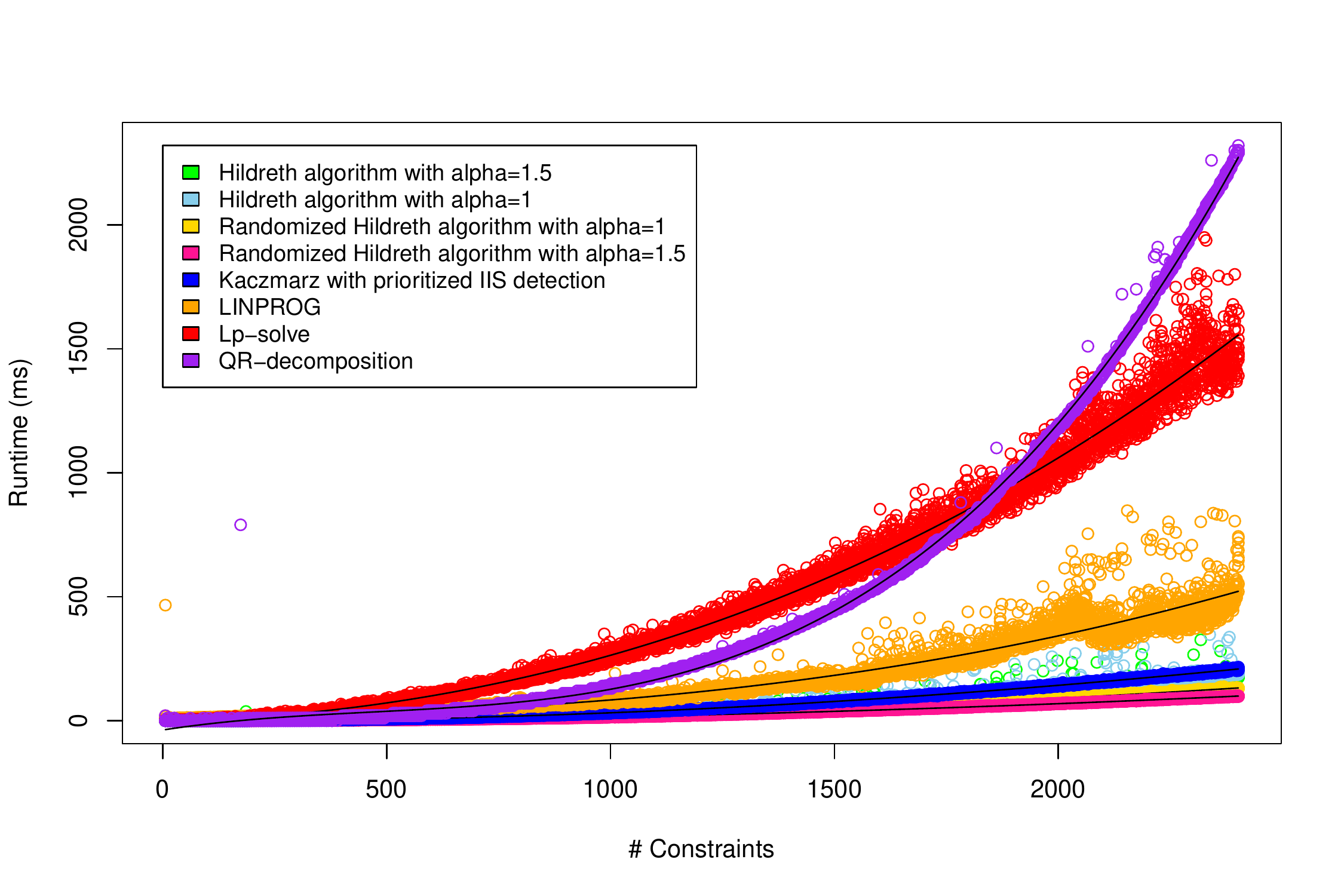}
\caption{Performance comparison of Hildreth algorithms, Randomized Hildreth algoritms, LINPROG, LP-Solve, and QR-decomposition}
\label{fig:lp-solve}
\end{figure}

\subsection{Methodology}

For all experiments we used the same hardware and test data generator, but instrumentalized the algorithms differently.
We used the following setup:
a desktop computer with Intel i5 3.3GHz processor and 64-bit Windows 7, running an Oracle Java virtual machine.
Layout specifications were randomly generated using the test data generator described in~\cite{Weber:High-Level-Constraints}.
For each experiment the same set of test data was used.
The specification size was varied from 4 to 2402 constraints, in increments of 4 constraints (2 new constraints for positioning and 2 new constraint for the preferred size of a new widget).
For each size, 10 different layouts were generated, resulting in a total of 6000 different layout specifications.
A tolerance of 0.01 was used for solving.

In the first experiment we investigated the convergence behavior of the algorithms.
We measured the number of sub-optimal solutions for each algorithm.
A solution is sub-optimal if the error of a constraint (the difference between the right-hand and left-hand side) is bigger than the given tolerance.

In the second experiment we measured the performance in terms of computation time (\(T\)) in milliseconds (ms), depending on the problem size measured in number of constraints~($c$).
The proposed algorithm was used to solve each of the problems of the test data set and the time was measured.
As a reference, all the generated specifications were also solved with an implementation of the cyclic Algorithm 2 (referred as Hildreth algorithm here), the randomized Algorithm 2 (referred as randomized Hildreth's algorithm here), Matlab's LINPROG solver~\cite{linprog:solver}, LP-Solve~\cite{lp:solve}, and the Kaczmarz with prioritized IIS detection~\cite{Gerald:Kaczmarz}.
LINPROG is widely known for its speed, and LP-Solve was previously used to solve UI layout problems~\cite{Weber:High-Level-Constraints}.
Additionally, we wanted to compare our algorithm with a direct method, so we also included the implementation of QR-decomposition in the Apache Commons Mathematics Library~\cite{Commons:Math} in the evaluation. Moreover, we test both versions of Algorithm 2 with two relaxation parameters $\alpha=1$ and $\alpha=1.5$.

\begin{table}[t]
\begin{center}
\begin{tabular}{rl}
  \toprule
 Symbol & Explanation \\
  \midrule
 \(\beta_0\) & Intercept of the regression model\\
 \(\beta_{1-3}\) & Estimated model parameters\\
 \(c\) & Number of constraints\\
  \(T\) & Measured time in milliseconds\\
   \(R^2\) & Coefficient of determination of the regression models\\
\bottomrule
\end{tabular}
\vspace{0.05in}
\caption{Symbol table}
\label{tab:symbs}
\end{center}
\end{table}

\newcommand{\sign}[1]{\textsuperscript{#1}}
\begin{table*}[t]
\begin{center}
\begin{tabular}{rlllll}
  \toprule
\multicolumn{1}{c}{\bf Strategy}                       &\multicolumn{1}{c}{\(\beta_0\)} &\multicolumn{1}{c}{\(\beta_1\)} & \multicolumn{1}{c}{\(\beta_2\)}  & \multicolumn{1}{c}{\(\beta_3\)} & \(R^2\)  \\
  \midrule
Hildreth algorithm with alpha=1.5 &\(\;\;\,1.799 \)\sign{***} & \(-1.096 \cdot 10^{-02}\)\sign{***} & \(\;\;\, 3.637 \cdot 10^{-05}\)\sign{***} & \(-4.369 \cdot 10^{-10}\) &\(0.9868\) \\
Hildreth algorithm with alpha=1 &\(\;\;\, 1.927 \)\sign{***} & \(-1.256\cdot 10^{-02}\)\sign{***} & \(\;\;\, 4.129 \cdot 10^{-05}\)\sign{***} & \( -1.177 \cdot 10^{-9}\)\sign{**}&\(0.9763\) \\
Randomized Hildreth algorithm with alpha=1 &\(\;\;\, 1.156 \)\sign{***} & \(-7.889 \cdot 10^{-03}\)\sign{***} & \(\;\;\, 2.737 \cdot 10^{-05}\)\sign{***} & \(-5.928 \cdot 10^{-10}\)\sign{***}&\(0.9994\) \\
Randomized Hildreth algorithm with alpha=1.5 &\(\;\;\, 1.202 \)\sign{***} & \(-7.331 \cdot 10^{-03}\)\sign{***} & \(\;\;\, 2.277 \cdot 10^{-05}\)\sign{***} & \(-1.138 \cdot 10^{-9}\)\sign{***}&\(0.9994\) \\
Kaczmarz with prioritized IIS detection &\(\;\;\, 1.035 \)\sign{***} & \( -1.112 \cdot 10^{-02}\)\sign{***} & \(\;\;\, 4.278 \cdot 10^{-05}\)\sign{***} & \(-9.176 \cdot 10^{-10}\)\sign{***}&\(0.9994\) \\
LINPROG &\(\;\;\, 18.29\)\sign{***} & \(\;\;\, 1.591 \cdot 10^{-04}\) & \(\;\;\, 4.934 \cdot 10^{-05}\)\sign{***} & \(\;\;\, 1.577 \cdot10^{-08}\)\sign{***} &\(0.9367\) \\
LP-Solve &\( -2.491\)\sign{***} & \(\;\;\, 3.924 \cdot 10^{-02}\)\sign{***} & \(\;\;\, 2.079 \cdot 10^{-04}\)\sign{***} & \(\;\;\, 1.904 \cdot 10^{-08}\)\sign{***} &\(0.9900\) \\
QR-Decomposition  &\( -37.70\)\sign{***} & \(\;\;\, 0.2802 \)\sign{***} & \(-4.009 \cdot 10^{-04}\)\sign{***} & \(\;\;\, 2.850 \cdot 10^{-07}\)\sign{***} & \(0.9989\) \\
\bottomrule
\multicolumn{6}{l}{Significance codes: \sign{***} \(p<0.001\)} \\
\end{tabular}
\vspace{0.05in}
\caption{Regression models for the different solving strategies}
\label{tab:reg1}
\end{center}
\end{table*}

\subsection{Results}

The first experiment tested the convergence behavior of the algorithms.
We found that all algorithms converge, which is expected since the algorithms were designed to find a solvable subproblem.

In the second experiment we investigated the performance behavior of the algorithms.
To identify the performance trend of the algorithms over $c$, we defined some regression models (linear, quadratic, log, cubic).
We found that the best-fitting model is the polynomial model
\[
T(c) = \beta_0 + \beta_1c + \beta_2c^2 + \beta_3c^3 + \epsilon,
\]
which gave us a good fit for the performance data.
Table~\ref{tab:symbs} explains the symbols used.
Key parameters of the models are depicted in Table~\ref{tab:reg1}; a graphical representation of the models can be found in Figure~\ref{fig:lp-solve}. Figure~\ref{fig:lp-solve} compares the aforementioned algorithms to LINPROG, LP-Solve, QR-decomposition, and Kaczmarz with prioritized IIS detection.
The Hildreth's algorithm with prioritized IIS detection (Algorithm 2) perform significantly better than LINPROG, LP-Solve and QR-decomposition, especially for bigger problems, and also outperforms the Kaczmarz with prioritized IIS detection.

\begin{figure}[t]
\centering
\includegraphics[width=\columnwidth]{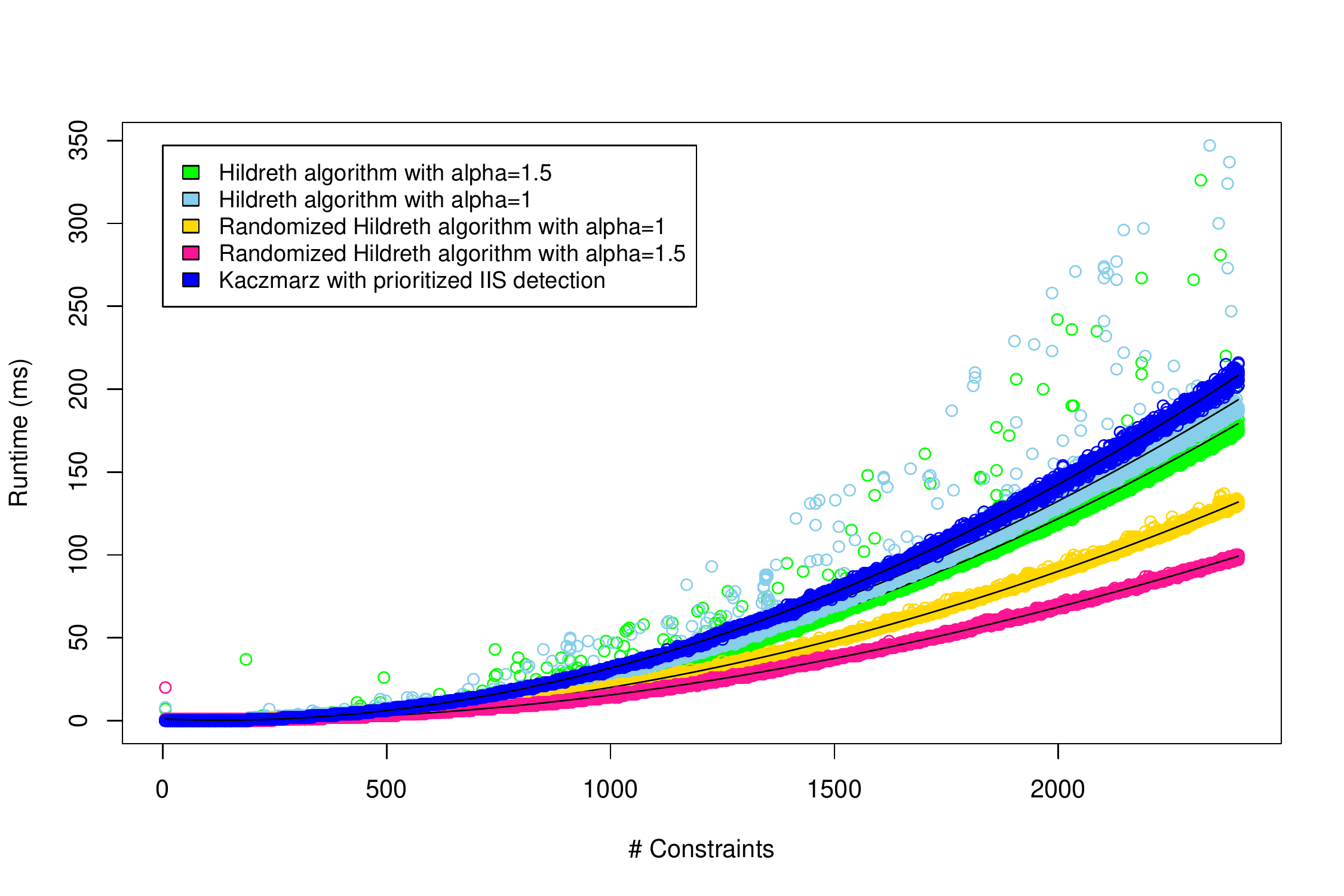}
\caption{Performance comparison of  Hildreth and randomized Hildreth algorithms}
\label{fig:random}
\end{figure}

Figure~\ref{fig:random} is an enlargement of Figure \ref{fig:lp-solve} with only the performance of the Hildreth's algorithm with prioritized IIS detection, and the Kaczmarz with prioritized IIS detection.
As the graphs show, randomized Algorithm 2 exhibits a better performance than the cyclic version. Moreover, over projection ($\alpha>1$) appears to converge faster as also observed in \cite{SV09:Randomized-Kaczmarz}.

%
%

\bibliographystyle{abbrv}
\bibliography{bibliography}
\end{document}